\documentclass[aos,seceqn,dvips]{arximspdf}
\usepackage{graphicx}

\doi{10.1214/10-AOS833}
\volume{38}
\issue{6}
\pubyear{2010}
\firstpage{3724}
\lastpage{3750}

\makeatletter

\newcommand{\eqref}[1]{(\ref{#1})}

\newtheorem{theorem}{Theorem}

\newtheorem{lemma}{Lemma}
\newproclaim{remark}{Remark}
\newtheorem{corollary}{Corollary}

\newcommand{\bbA}{\mathbf{A}}
\newcommand{\bbB}{\mathbf{B}}
\newcommand{\bbc}{\mathbf{c}}
\newcommand{\bbD}{\mathbf{D}}
\newcommand{\bbF}{\mathbf{F}}
\newcommand{\bbH}{\mathbf{H}}
\newcommand{\bbI}{\mathbf{I}}
\newcommand{\bbS}{\mathbf{S}}
\newcommand{\bbT}{\mathbf{T}}
\newcommand{\bbW}{\mathbf{W}}
\newcommand{\bbX}{\mathbf{X}}
\newcommand{\bbY}{\mathbf{Y}}

\makeatother

\begin{document}
\begin{frontmatter}

\title{Nonparametric estimate of spectral density functions of sample covariance matrices: A~first~step}
\runtitle{Nonparametric estimate of density functions}

\begin{aug}
\author[A]{\fnms{Bing-Yi} \snm{Jing}\thanksref{t1}\ead[label=e1]{majing@ust.hk}},
\author[B]{\fnms{Guangming} \snm{Pan}\thanksref{t2}\ead[label=e2]{gmpan@ntu.edu.sg}},
\author[A]{\fnms{Qi-Man}~\snm{Shao}\thanksref{t3}\ead[label=e3]{maqmshao@ust.hk}}\break
\and
\author[C]{\fnms{Wang}~\snm{Zhou}\thanksref{t4}\corref{}\ead[label=e4]{stazw@nus.edu.sg}}
\runauthor{Jing, Pan, Shao and Zhou}
\thankstext{t1}{Supported in part by Hong Kong RGC Grants
HKUST6011/07P and HKUST6015/08P.}
\thankstext{t2}{Supported in part by Grant M58110052 at the Nanyang
Technological University.}
\thankstext{t3}{Supported in part by Hong Kong RGC CERG 602608.}
\thankstext{t4}{Supported in part by Grant R-155-000-083-112 at the
National University of Singapore.}
\affiliation{Hong Kong University of Science and Technology, Nanyang
Technological University,
Hong Kong University of Science and Technology and~National~University~of~Singapore}

\address[A]{B.-Y. Jing\\ 
Q.-M. Shao\\
Department of Mathematics\\
Hong Kong University\\
\quad of Science and Technology\\
Clear Water Bay\\
Kowloon\\
Hong Kong\\
\printead{e1}\\
\phantom{E=mail:}\printead*{e3}}

\address[B]{G. Pan\\
Division of Mathematical Sciences\\
School of Physical\\
\quad and Mathematical Sciences\\
Nanyang Technological University\\
Singapore 637371\\
\printead{e2}}

\address[C]{W. Zhou\\
Department of Statistics\\
\quad and Applied Probability\\
National University of Singapore\\
Singapore 117546\\
\printead{e4}}
\end{aug}

\received{\smonth{1} \syear{2010}}
\revised{\smonth{5} \syear{2010}}

%
\begin{abstract}
The density function of the limiting spectral distribution of general
sample covariance matrices is usually
unknown. We propose to use kernel estimators which are proved to be
consistent. A~simulation study is also
conducted to show the performance of the estimators.
\end{abstract}

%

\end{frontmatter}

\section{Introduction}

Suppose that $X_{ij}$ are independent and identically distributed
(i.i.d.) real random variables. Let $\bbX_n=(X_{ij})_{p\times n}$ and
$\bbT_n$ be a $p\times p$ nonrandom Hermitian nonnegative definite
matrix. Consider the random matrices
\[
\bbA_n=\frac{1}{n}\bbT_n^{1/2}\bbX_n\bbX_n^T\bbT_n^{1/2}.
\]
When $EX_{11}=0$ and $EX_{11}^2=1$, $\bbA_n$ can be
viewed as a sample covariance matrix drawn from the population with
covariance matrix $\bbT_n$. Moreover, if $\bbT_n$ is another sample
covariance matrix, independent of $\bbX_n$, then $\bbA_n$ is a Wishart
matrix.

Sample covariance matrices are of paramount importance in multivariate
analysis. For example, in principal component analysis, we need to
estimate eigenvalues of sample covariance matrices in order to obtain
an interpretable low-dimensional data representation. The matrices
consisting of contemporary data are usually large, with the number of
variables proportional to the sample size. In this setting, fruitful
results have accumulated since the celebrated Marcenko and Pastur
law~\cite{MP} was discovered; see the latest monograph of Bai and
Silverstein~\cite{bai06} for more details.

The basic limit theorem regarding $\bbA_n$ concerns its empirical
spectral distribution $F^{\bbA_n}$. Here, for any matrix $\bbA$ with
real eigenvalues, the empirical spectral distribution $F^{\bbA}$ is
given by
\[
F^{\bbA}(x)=\frac{1}{p}\sum_{k=1}^pI(\lambda_k\leq x),
\]
where $\lambda_k$, $k=1,\ldots,p$, denote the eigenvalues of $\bbA$.

Suppose the ratio of the dimension to the sample size $c_n=p/n$ tends
to $c$ as $n\to\infty$. When $\bbT_n$ becomes the identity matrix,
$F^{\bbA_n}$ tends to the so-called Marcenko and Pastur law with the
density function
\[
f_c(x)=\cases{
(2\pi cx)^{-1}\sqrt{(b-x)(x-a)}, &\quad \vspace*{2pt}$a\leq x\leq b$,\cr
0, &\quad otherwise.
}
\]
It has point mass $1-c^{-1}$ at the origin if $c>1$, where $a=(1-\sqrt
c)^2$ and $b=(1+\sqrt c)^2$ (see Bai and Silverstein~\cite{bai06}).

In the literature, it is also common to study
\[
\bbB_n=\frac{1}{n}\bbX_n^T\bbT_n\bbX_n
\]
since the eigenvalues of $\bbA_n$ and $\bbB_n$ differ by $|n-p|$ zero
eigenvalues. Thus,
%
\begin{equation} \label{relation}
F^{\bbB_n}(x)=\biggl(1-\frac{p}{n}\biggr)I\bigl(x\in
[0,\infty)\bigr)+\frac{p}{n}F^{\bbA_n}(x).
\end{equation}
When $F^{\bbT_n}$ converges weakly to a nonrandom distribution $H$,
Marcenko and Pastur~\cite{MP}, Yin~\cite{y1} and Silverstein~\cite{s3}
proved that, with probability one, $F^{\bbB_n}(x)$ converges in
distribution to a nonrandom distribution function
$\underline{F}_{c,H}(x)$
whose Stieltjes transform $\underline m(z)=m_{\underline{F}_{c,H}}(z)$
is, for each $z\in\mathcal{C}^+=\{z\in\mathcal{C}\dvtx \Im z>0\}$, the
unique solution to
the equation
%
\begin{equation}\label{a3}
\underline{m}=-\biggl(z-c\int\frac{t\,dH(t)}{1+t\underline{m}}\biggr)^{-1}.
\end{equation}
Here, the Stieltjes transform $m_F(z)$ for any probability
distribution function $F(x)$ is defined by
%
\begin{equation}\label{b5}
m_F(z)=\int\frac{1}{x-z}\,dF(x),\qquad z\in\mathcal{C}^+.
\end{equation}
Therefore,
from~(\ref{relation}), we have
%
\begin{equation}\label{b4}
\underline{F}_{c,H}(x)=(1-c)I\bigl(x\in[0,\infty)\bigr)+cF_{c,H}(x),
\end{equation}
where $F_{c,H}(x)$ is the limit of $F^{\bbA_n}(x)$. As a consequence of
this fact, we have
%
\begin{equation}\label{b3}
\underline{m}(z)=-\frac{1-c}{z}+cm(z).
\end{equation}
Moreover, $\underline m(z)$ has an inverse,
%
\begin{equation}\label{a4}
z(\underline{m})=-\frac{1}{\underline{m}}+c_n\int\frac{t\,dH(t)}{1+t\underline{m}}.
\end{equation}
Relying on this inverse, Silverstein and Choi~\cite{s1} carried out a
remarkable analysis of the analytic behavior of
$\underline{F}_{c,H}(x)$.

When $\bbT_n$ becomes the identity matrix, there is an explicit
solution to~(\ref{a3}). In this case, from \eqref{relation}, we see
that the density function of $\underline{F}_{c,H}(x)$ is
\[
\underline{f}_{c,I}(x)=(1-c)I(c<1)\delta_0+cf_c(x),
\]
where $\delta_0$ is the point mass at $0$. Unfortunately, there is no
explicit solution to~(\ref{a3}) for general $\bbT_n$. Although we can
use $F^{\bbA_n}(x)$ to estimate $F_{c,H}(x)$, we cannot make any
statistical inference on $F_{c,H}(x)$ because there is, as\vspace*{1pt} far as we
know, no central limit theorem concerning
$(F^{\bbA_n}(x)-F_{c,H}(x))$. Actually, it is argued in Bai and
Silverstein~\cite{bai06} that the process
$n(F^{\bbA_n}(x)-F_{c,H}(x))$, $x\in(-\infty, \infty),$ does
not converge to a nontrivial process in any metric space. This makes us
want to pursue other ways of understanding the limiting spectral
distribution $F_{c,H}(x)$.

This paper is part of a program to estimate the density function
$f_{c,H}(x)$ of the limiting spectral distribution $F_{c,H}(x)$ of
sample covariance matrices $\bbA_n$ by kernel estimators. In this
paper, we will prove the consistency of those estimators as a first
step.

\section{Methodology and main results}

Suppose that the observations $X_1,\ldots,\break X_n$ are i.i.d. random
variables with an unknown density function $f(x)$ and $F_n(x)$ is the
empirical distribution function determined by the sample. A~popular
nonparametric estimate of $f(x)$ is then
%
\begin{equation}\label{a24}
\hat{f}_n(x)=\frac{1}{nh}\sum_{j=1}^nK\biggl(\frac
{x-X_j}{h}\biggr)=\frac{1}{h}\int
K\biggl(\frac{x-y}{h}\biggr)\,dF_n(y),
\end{equation}
where the function $K(y)$ is a Borel function and $h=h(n)$ is the
bandwidth which tends to $0$ as $n\to\infty$. Obviously,
$\hat{f}_n(x)$ is again a probability density function and, moreover,
it inherits some smooth properties of $K(x)$, provided the kernel is
taken as a probability density function. Under some\vspace*{1pt} regularity
conditions on the kernel, it is well known that
$\hat{f}_n(x)\rightarrow f(x)$ in some sense (with probability one, or
in probability). There is a huge body of literature regarding this kind
of estimate. For example, one may refer to Rosenblatt~\cite{R1}, Parzen
\cite{p1}, Hall~\cite{peter} or the book by Silverman~\cite{bw}.

Informed by~(\ref{a24}), we propose the following estimator $f_n(x)$ of
$f_{c,H}(x)$:
%
\begin{equation}\label{a1}
f_n(x)=\frac{1}{ph}\sum_{i=1}^pK\biggl(\frac{x-\mu_i}{h}\biggr)=\frac
{1}{h}\int
K\biggl(\frac{x-y}{h}\biggr)\,dF^{\bbA_n}(y),
\end{equation}
where $\mu_i$, $i=1,\dots,p$, are eigenvalues of $\bbA_n$. It turns out
that $f_n(x)$ is a consistent estimator of $f_{c,H}(x)$ under some
regularity conditions.

Suppose that the kernel function $K(x)$ satisfies
%
\begin{equation}\label{a25}
\sup_{-\infty< x<\infty} |K(x)|<\infty,\qquad\lim_{|x|\rightarrow\infty}|xK(x)|=0\
\end{equation}
and
%
\begin{equation} \label{a26}
\int K(x)\,dx=1,\qquad\int|K'(x)|\,dx<\infty.
\end{equation}

\begin{theorem}\label{theo1}
Suppose that $K(x)$ satisfies
(\ref{a25}) and~(\ref{a26}). Let $h = h(n)$ be a sequence of positive
constants satisfying
%
\begin{equation}\label{band}
\lim_{n\to\infty}nh^{5/2}=\infty, \qquad\lim_{n\to\infty}h=0.
\end{equation}
Moreover, suppose that all $X_{ij}$
are i.i.d. with $EX_{11}=0$, $\operatorname{Var}(X_{11})=1$ and $EX_{11}^{16}<\infty$.
Also, assume that
$c_n\rightarrow c\in(0,1)$.
Let $\bbT_n$ be a $p\times p$ nonrandom symmetric positive definite
matrix with spectral norm bounded above by a positive constant such
that $H_n=F^{\bbT_n}$ converges weakly to a nonrandom distribution $H$.
In addition, suppose that $F_{c,H}(x)$ has a compact support $[a,b]$
with $a>0$. Then,
\[
f_n(x)\longrightarrow f_{c,H}(x)\qquad\mbox{in probability uniformly in } x\in[a,b].
\]

\end{theorem}

\begin{remark} \label{remark}
We conjecture that the condition $EX_{11}^{16}$ can be reduced to
$EX_{11}^4<\infty$.
\end{remark}

When $\bbT_n$ is the identity matrix, we have a slightly better result.
\begin{theorem} \label{theo2} Suppose that $K(x)$ satisfies
(\ref{a25}) and~(\ref{a26}). Let $h=h(n)$ be a sequence of positive
constants satisfying
%
\begin{equation}\label{band1}
\lim_{n\rightarrow\infty} nh^2=\infty,\qquad
\lim_{n\rightarrow\infty} h=0.
\end{equation}
Moreover,
suppose that all $X_{ij}$ are i.i.d. with $EX_{11}=0$, $\operatorname{Var}(X_{11})=1$
and $EX_{11}^{12}<\infty$. Also, assume that
$c_n\rightarrow c\in(0,1)$.
Denote the support of the MP law by $[a,b]$. Let $\bbT_n=\bbI$. Then,
\[
\sup_{x\in[a,b]}|f_n(x)- f_{c}(x)|\longrightarrow0\qquad\mbox{in probability}.
\]
\end{theorem}

Theorem~\ref{theo1} also gives the estimate of $F_{c,H}(x),$ as
below.\vadjust{\goodbreak}

\begin{corollary} \label{cor1}
Under the assumptions of Theorem
\ref{theo1}, correspondingly,
%
\begin{equation}\label{a31}
F_n(x)\rightarrow F_{c,H}(x)\qquad\mbox{in probability},
\end{equation}
where
%
\begin{equation}\label{c6}
F_n(x)=\int^x_{-\infty}f_n(t)\,dt.
\end{equation}
\end{corollary}

Corollary~\ref{cor1} and the Helly--Bray lemma ensure that we have the
following.

\begin{corollary} \label{cor2}
Under the assumptions of Theorem
\ref{theo1}, if $g(x)$ is a continuous bounded function, then
%
\begin{equation}\label{a31}
\int g(x)\,dF_n(x)\rightarrow\int g(x)\,dF_{c,H}(x)\qquad\mbox{in probability}.
\end{equation}
\end{corollary}

In order to prove consistency of the nonparametric estimates, we need
to develop a convergence rate for $F^{\bbA_n}$. When $\bbT_n=\bbI$, Bai
\cite{b3} developed a Berry--Esseen-type inequality and investigated
the convergence rate of $EF^{\bbA_n}$. Later, G\"otze and Tikhomirov
\cite{g2} improved the Berry--Esseen-type inequality and obtained a
better convergence rate. For general $\bbT_n$, we establish the
following convergence rate.

\begin{theorem}\label{theo3}
Under the assumptions of Theorem~\ref{theo1},
%
\begin{equation}\label{b34}
\sup_{x}|EF^{\bbA_n}(x)-F_{c_n,H_n}(x)|=O\biggl(\frac{1}{n^{2/5}}\biggr)
\end{equation}
and
%
\begin{equation}\label{b35}
E\sup_{x}|F^{\bbA_n}(x)-F_{c_n,H_n}(x)|=O\biggl(\frac{1}{n^{2/5}}\biggr).
\end{equation}
\end{theorem}

\begin{remark}
Under the fourth moment condition, that is, $EX_{11}^4<\infty$, we
conjecture that the above rate $O(n^{-2/5})$ could be improved to
$O(n^{-1}\sqrt{\log n})$.
\end{remark}

\section{Applications}

Let us demonstrate some applications of Theorems~\ref{theo1},~\ref{theo2} and their corollaries. Since $F_{c,H}(x)$ does not have an
explicit expression (except for some special cases), we may now use
$F_n(x)$ to estimate it, by Corollary~\ref{cor1}. More importantly,
$F_n(x)$ has some smoothness properties, which $F^{\bbA_n}$ does not
have.

We first consider an example in wireless communication. Consider a
synchronous CDMA system with $n$ users and processing gain $p$. The
discrete-time model for the received signal $\bbY$ is given by
%
\begin{equation}\label{1}
\bbY=\sum_{k=1}^nx_k\mathbf{h}_k+\bbW,
\end{equation}
where $x_i\in\mathcal{R}$ and $\mathbf{h}_k\in\mathcal{R}^p$ are,
respectively, the transmitted symbol and the signature spreading
sequence of user $k$, and $\bbW$ is the Gaussian noise with zero mean
and covariance matrix $\sigma^2\bbI$. Assume that the transmitted
symbols of different users are independent, with $Ex_k=0$ and
$E|x_k|^2=p_k$. This model is slightly more general than that in
\cite{Tse99}, where all of the users' powers $p_k$ are assumed to be
the same.

Following~\cite{Tse99}, consider the demodulation of user $1$ and use
the signal-to-interference ratio (SIR) as the performance measure of
linear receivers. The SIR of user $1$ is defined by (see~\cite{Tse99})
\[
\beta_1=\frac{(\bbc_1^T\mathbf{h}_1)^2p_1}{\bbc_1^T\bbc_1\sigma
^2+\sum_{k=2}^K(\bbc_1^T\mathbf{h}_k)^2p_k}.
\]
The minimum mean square error (MMSE) receiver minimizes the mean square
error as well as maximizes the SIR for all users (see~\cite{Tse99}).
The SIR of user $1$ is given by
\[
\beta_1^{\mathrm{MMSE}}=p_1\mathbf{h}_1^T(\bbH_1\bbD_1\bbH_1^T+\sigma^2\bbI
)^{-1}\mathbf{h}_1,
\]
where
\[
\bbD_1=\operatorname{diag}(p_2,\ldots,p_n),\qquad \bbH_1=(\mathbf{h}_2,\ldots,\mathbf{h}_n).
\]
Assume that the $\mathbf{h}_k'$ are i.i.d. random vectors, each consisting of
i.i.d. random variables with appropriate moments. Moreover, suppose
that $p/n\rightarrow c>0$ and $F^{\bbD_1}(x)\rightarrow H(x)$. Then, by
Lemma 2.7 in~\cite{b4} and the Helly--Bray lemma, it is not difficult
to check that
\[
\beta_1^{\mathrm{MMSE}}-p_1\int\frac{1}{x+\sigma^2}\,dF_{c,H}(x)\stackrel
{i.p.}\longrightarrow
0.
\]
To judge the performance of different receivers, we may then compare
the value of $\int\frac{1}{x+\sigma^2}\,dF_{c,H}(x)$ with the limiting
SIR of the other linear receiver.
However, the awkward fact is that we usually do not have an explicit
expression for $F_{c,H}(x)$. Thus, we may use the kernel estimate
$\int\frac{1}{x+\sigma^2}\,dF_n(x)$ to estimate
$\int\frac{1}{x+\sigma^2}\,dF_{c,H}(x)$, by Corollary~\ref{cor2}.

A second application: we may use $f_{n}(x)$ to infer, in some way, some
statistical properties of the population covariance matrix $\bbT_n$.
Specifically speaking, by~(\ref{b5}), we may evaluate the Stieltjes
transform of the kernel estimator $f_n(x),$
%
\begin{equation}\label{b6}
m_{f_n}(z)=\int\frac{1}{x-z}f_n(x)\,dx,\qquad z\in\mathcal{C}^+.
\end{equation}
We may then obtain $\underline{m}_{f_n}(z),$ by~(\ref{b3}). On the
other hand, we conclude from~(\ref{a4}) that
%
\begin{equation}\label{c1}
\frac{\underline{m}(z)(c-1-z\underline{m}(z))}{c}=\int\frac
{dH(t)}{t+1/\underline{m}(z)}.
\end{equation}
Note that $\underline{m}(z)$ has a positive imaginary part. Therefore,
with notation $z_1=-1/\underline{m}(z)$ and
$s(z_1)=\frac{\underline{m}(z)(c-1-z\underline{m}(z))}{c}$, we can
rewrite~(\ref{c1}) as
%
\begin{equation}\label{c3}
s(z_1)=\int\frac{dH(t)}{t-z_1},\qquad z_1\in\mathcal{C}^+.
\end{equation}
Consequently, in view of the inversion formula
%
\begin{equation}
\label{c2} F\{[a,b]\}=\frac{1}{\pi}\lim_{v\rightarrow
0}\int^b_a
\Im m_F(u+iv)\,du,
\end{equation}
we may recover $H(t)$ from $s(z_1)$ as given in~(\ref{c3}). However,
$s(z_1)$ can be estimated by the resulting kernel estimate
%
\begin{equation}
\label{c5} \frac{\underline
{m}_{f_n}(z)(c-1-z\underline{m}_{f_n}(z))}{c}.
\end{equation}
Once $H(t)$ is estimated, we may further estimate the functions of the
population covariance matrix $\bbT_n$, such as $\frac{1}{n}\operatorname{tr}\bbT_n^2$.
Indeed, by the Helly--Bray lemma, we have
\[
\frac{1}{n}\operatorname{tr}\bbT_n^2=\int t^2\,dH_n(t)\stackrel{D}\longrightarrow
\int t^2\,dH(t).
\]
Thus, we may construct an estimator for $\frac{1}{n}\operatorname{tr}\bbT_n^2$ based
on the resulting kernel estimate~(\ref{c5}). We conjecture that the
estimators of $H(t)$ and the corresponding functions like
$\frac{1}{n}\operatorname{tr}\bbT_n^2$, obtained by the above method, are also
consistent. A~rigorous argument is currently being pursued.

\section{Simulation study}
In this section, we perform a simulation study to investigate the
behavior of the kernel density estimators of the Marcenko and Pastur
law. We consider two different populations, exponential and binomial
distributions. From each population, we generate two samples with sizes
$50\times200$ and $800\times3200,$ respectively. We can therefore
form two random matrices, $(X_{ij})_{50,200}$ and
$(X_{ij})_{800,3200}$. The kernel is selected as
\[
K(x)=(2\pi)^{-1/2}e^{-x^2/2},
\]
which is the standard normal density function. The bandwidth is chosen
as $h=0.5n^{-1/3}$ ($n=200, 3200$).

For $(X_{ij})_{50,200}$, the kernel density estimator is
\[
\frac1{50\times200^{-2/5}}\sum_{i=1}^{50}
K\bigl((x-\mu_i)/200^{-2/5}\bigr),
\]
where $\mu_i, i=1,\ldots,50$, are eigenvalues of
$200^{-1}(X_{ij})_{50,200}(X_{ij})_{50,200}^T$. This curve is drawn by
dot-dash lines in the first two pictures.

For $(X_{ij})_{800,3200}$, the kernel density estimator is
\[
\frac1{800\times3200^{-2/5}}\sum_{i=1}^{800}
K\bigl((x-\mu_i)/3200^{-2/5}\bigr),
\]
where $\mu_i$, $i=1,\ldots,800$, are eigenvalues of
$3200^{-1}(X_{ij})_{800,3200}(X_{ij})_{800,3200}^T$. This curve is
drawn by dashed lines in the first two pictures.

\begin{figure}

\includegraphics{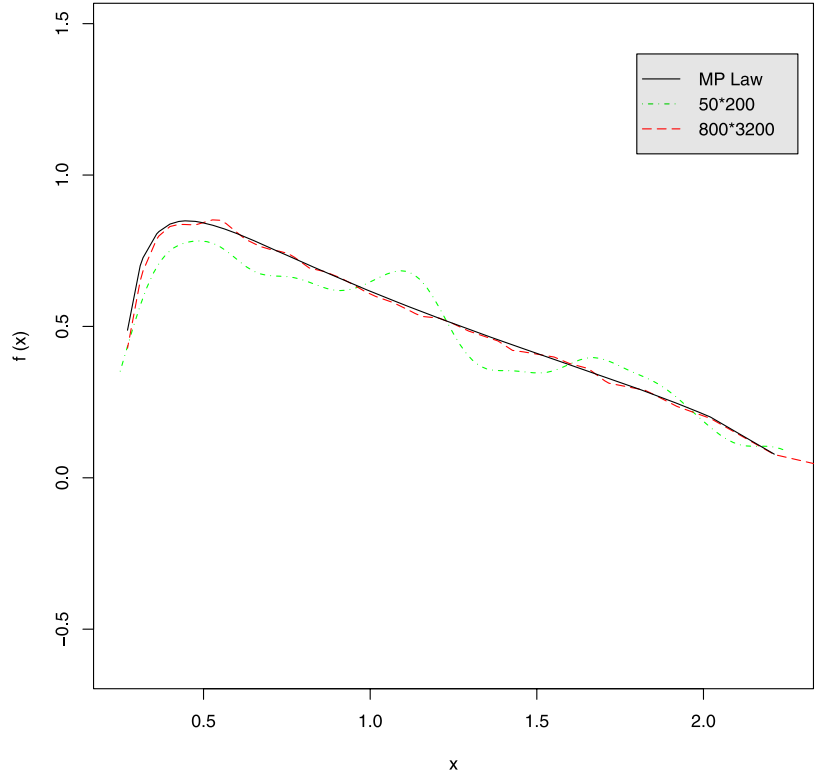}

\caption{Spectral density curves for sample covariance matrices
$n^{-1}(X_{ij})_{p\times n}(X_{ij})_{p\times n}^T$, $X_{ij}\sim$ exponential distribution.}\label{fig1}
\end{figure}

\begin{figure}

\includegraphics{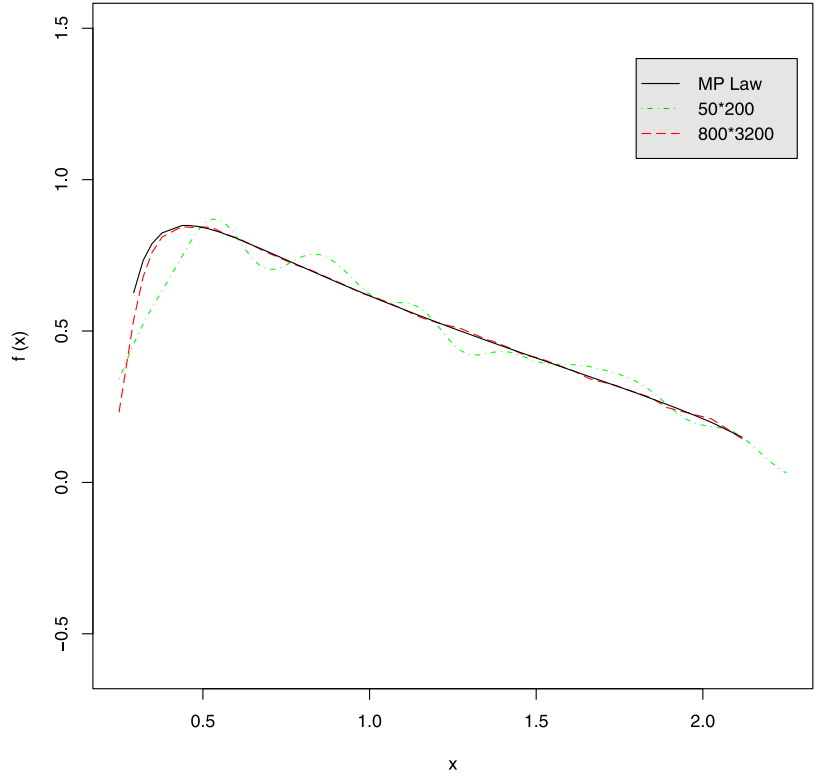}

\caption{Spectral density curves for sample covariance matrices
$n^{-1}(X_{ij})_{p\times n}(X_{ij})_{p\times n}^T$,
  $X_{ij}\sim$ binomial distribution.}\label{fig2}
\end{figure}

The density function of the Marcenko and Pastur law is drawn by solid
lines in the first two pictures. Here, in Figure~\ref{fig1}, the distribution is
%
\begin{equation} \label{exp}
F(x)=e^{-(x+1)},\qquad x\geq-1.
\end{equation}
In Figure~\ref{fig2}, the distribution is
%
\begin{equation} \label{bion}
P(X=-1)=1/2,\qquad P(X=1)=1/2.
\end{equation}
From the two figures, we see that the estimated curves fit the
Marcenko and Pastur law very well. As $n$ becomes large, the estimated
curves become closer to the Marcenko and Pastur law.

\begin{figure}

\includegraphics{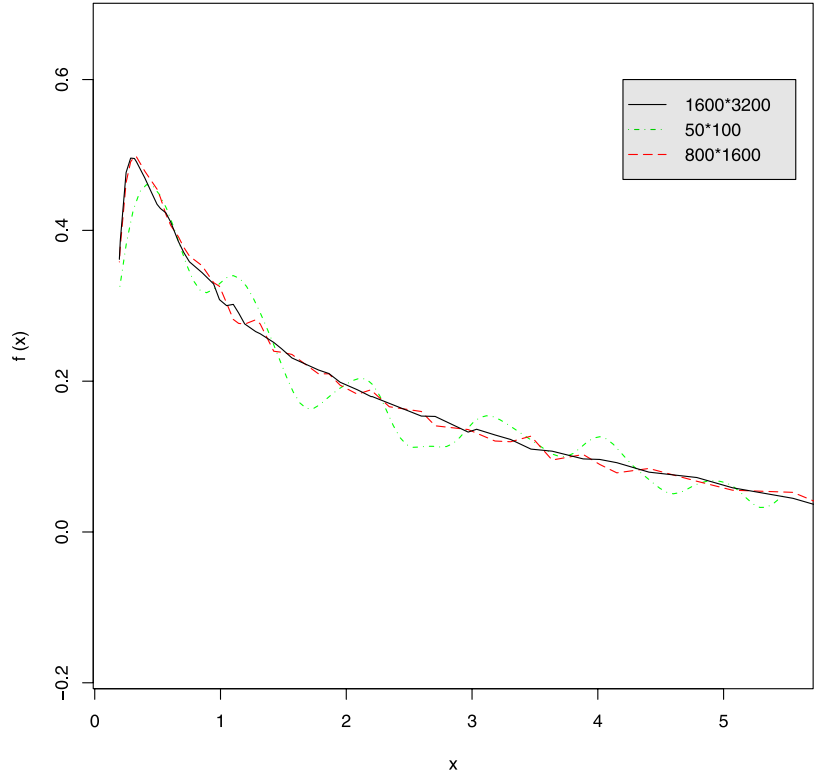}

\caption{Spectral density curves for sample covariance matrices
 $n^{-1}\bbT_n^{1/2}(X_{ij})_{p\times n}(X_{ij})_{p\times n}^T\bbT_n^{1/2}$,
 $X_{ij}\sim$ exponential distribution,
 $\bbT_n=\allowbreak n^{-1}(Y_{ij})_{p\times n}(Y_{ij})_{p\times n}^T$,
  $Y_{ij}\sim$ binomial distribution.}\label{fig3}\vspace*{12pt}
\end{figure}

Finally, we consider the estimated density curves based on the
following three matrices:
\begin{eqnarray*}
\bbA_{200}&=&\tfrac{1}{200}\bbT_{200}^{1/2}\bbX_{50\times
200}\bbX_{50\times200}^T\bbT_{200}^{1/2},\\
\bbA_{3200}&=&\tfrac{1}{3200}\bbT_{800}^{1/2}\bbX_{800\times
3200}\bbX_{800\times3200}^T\bbT_{3200}^{1/2},\\
\bbA_{6400}&=&\tfrac{1}{6400}\bbT_{6400}^{1/2}\bbX_{1600\times
6400}\bbX_{1600\times6400}^T\bbT_{6400}^{1/2},
\end{eqnarray*}
where $\bbX_{p\times4p}$, $p=50, 800, 1600$, are $p\times4p$ matrices
whose elements are i.i.d. random variables with
distribution \eqref{exp}, and $\bbT_{n}=\frac{1}{4p}\bbY_{p\times
4p}\bbY_{p\times4p}^T $. Here, $\bbY_{p\times4p}$ is a $p \times4p$
matrix consisting of i.i.d. random variables whose distributions are
given by \eqref{bion}. $\bbT_n$ and $\bbX_{p\times4p}$ are
independent. The kernel function is the same as before. The bandwidths
corresponding to the three matrices are $0.5\times(4p)^{-1/3}$. In
Figure~\ref{fig3}, we present three estimated curves. The dot-dash line is based
on $\bbA_{200}$, the dashed line on $\bbA_{3200}$ and the solid line
on $\bbA_{6400}$. Although, in this case, we do not know its exact
formula, we can predict the limiting spectral density function from
Figure~\ref{fig3}.

In order to show that the above conclusion is reliable, we choose ten
points throughout the range and calculate the mean square errors (MSEs)
for the kernel density estimator at the selected ten points, based on
$500$ matrices,\looseness=1
\[
\operatorname{MSE}(x)=500^{-1}\sum_{i=1}^{500}\bigl(f_n^{(i)}(x)-f_c(x)\bigr)^2,
\]\looseness=0
where $f_n^{(i)}(x)$ is the kernel density estimator at $x$ based on
the $i$th matrix. If the limiting distribution is unknown as in the
case $\bbA_{200}$, we use the averaged spectral density
\[
\bar f_c(x)=500^{-1}\sum_{i=1}^{500}f_n^{(i)}(x).
\]
So, in this case,
\[
\operatorname{MSE}(x)=500^{-1}\sum_{i=1}^{500}\bigl(f_n^{(i)}(x)-\bar f_c(x)\bigr)^2.
\]
The numerical results for the three different matrices considered in
this section are presented in Tables~\ref{tab1},~\ref{tab2} and~\ref{tab3}. The notation
``e$-$\textit{j}'' in these tables means multiplication by $10^{-j}$. The
MSEs are uniformly small. As $n$ becomes large, the MSEs become
smaller. This supports the conclusion that our proposed kernel spectral
density curve is consistent.

\begin{table}\tablewidth=280pt
\caption{MSE of spectral density curves for sample covariance matrices
 $n^{-1}(X_{ij})_{p\times n}(X_{ij})_{p\times n}^T$, $X_{ij}\sim$ exponential
 distribution}\label{tab1}
\begin{tabular*}{280pt}{@{\extracolsep{\fill}}lccccc@{}}
\hline
\multicolumn{1}{@{}r}{$\bolds{x=}$} &  \textbf{0.30}  & \textbf{0.511}   & \textbf{0.722}   & \textbf{0.933}   & \textbf{1.144} \\
\hline
\phantom{0}$50\times 200$ &9.89e$-$2\phantom{0}  &3.21e$-$2  &3.18e$-$2  &3.25e$-$2  &3.56e$-$2  \\
$800\times 3200$
&3.84e$-$03  &7.44e$-$5  &7.28e$-$5  &7.67e$-$5  &7.34e$-$5    \\ [6pt]
\multicolumn{1}{@{}r}{$\bolds{x=}$} & \textbf{1.356}   & \textbf{1.567} &\textbf{1.778}  &\textbf{1.989} &\textbf{2.20} \\
\hline
\phantom{0}$50\times 200$&3.79e$-$2\phantom{0} &3.18e$-$2 &3.73e$-$2 &2.76e$-$2 &3.63e$-$2 \\
$800\times 3200$
&7.67e$-$5\phantom{0} &7.23e$-$5 &6.88e$-$5 &6.60e$-$5 &6.74e$-$5\\
\hline
\end{tabular*}
\end{table}

\begin{table}[b]
\tablewidth=280pt
\caption{MSE of spectral density curves for sample covariance matrices
$n^{-1}(X_{ij})_{p\times n}(X_{ij})_{p\times n}^T$, $X_{ij}\sim$ binomial
distribution}\label{tab2}
\begin{tabular*}{280pt}{@{\extracolsep{\fill}}lccccc@{}}
\hline
\multicolumn{1}{@{}r}{$\bolds{x=}$} &  \textbf{0.30}  & \textbf{0.511}   & \textbf{0.722}   & \textbf{0.933}   & \textbf{1.144} \\
\hline
\phantom{0}$50\times 200$ &3.23e$-$1  &3.14e$-$2  &2.38e$-$2  &2.76e$-$2  &2.86e$-$2   \\
$800\times 3200$
&\phantom{0}5.13e$-$03  &8.01e$-$5  &6.05e$-$5  &7.30e$-$5  &6.53e$-$5    \\ [6pt]
\multicolumn{1}{@{}r}{$\bolds{x=}$} & \textbf{1.356}   & \textbf{1.567} &\textbf{1.778}  &\textbf{1.989} &\textbf{2.20} \\
\hline
\phantom{0}$50\times 200$&2.70e$-$2 &2.44e$-$2 &2.42e$-$2 &2.40e$-$2 &1.69e$-$2 \\
$800\times 3200$
&6.28e$-$5 &7.65e$-$5 &6.14e$-$5 &6.68e$-$5 &1.13e$-$4\\
\hline
\end{tabular*}
\end{table}

\begin{table}\tablewidth=280pt
\caption{MSE of spectral density curves for sample covariance matrices
$n^{-1}\bbT_n^{1/2}(X_{ij})_{p\times n}(X_{ij})_{p\times n}^T\bbT_n^{1/2}$, $X_{ij}\sim$ exponential distribution
$\bbT_n=n^{-1}(Y_{ij})_{p\times n}(Y_{ij})_{p\times n}^T$, $Y_{ij}\sim$ binomial distribution
$n^{-1}(X_{ij})_{p\times n}(X_{ij})_{p\times n}^T$, $X_{ij}\sim$ binomial
distribution}\label{tab3}
\begin{tabular*}{280pt}{@{\extracolsep{\fill}}lccccc@{}}
\hline
\multicolumn{1}{@{}r}{$\bolds{x=}$} &  \textbf{0.30}  & \textbf{0.511}   & \textbf{0.722}   & \textbf{0.933}   & \textbf{1.144} \\
\hline
\phantom{00}$50\times 100$ &1.20e$-$2  &8.71e$-$3  &8.58e$-$3  &7.90e$-$3  &8.77e$-$3   \\
\phantom{0}$800\times 1600$
     &\phantom{0}6.25e$-$05  &4.00e$-$5  &3.51e$-$5  &3.19e$-$5  &2.71e$-$5    \\
$1600\times 3200$
     & 2.98e$-$5 &1.83e$-$5 &1.44e$-$5 &1.39e$-$5 &1.53e$-$5      \\ [6pt]
\multicolumn{1}{@{}r}{$\bolds{x=}$} & \textbf{1.356}   & \textbf{1.567} &\textbf{1.778}  &\textbf{1.989} &\textbf{2.20} \\
\hline
\phantom{00}$50\times 200$&7.91e$-$3 &8.07e$-$3 &8.34e$-$3 &7.54e$-$3 &7.17e$-$3 \\
\phantom{0}$800\times 3200$
&3.04e$-$5 &3.10e$-$5 &2.98e$-$5 &2.89e$-$5 &2.66e$-$5\\
$1600\times 3200$ & 1.19e$-$5 &1.19e$-$5 &1.36e$-$5 &1.29e$-$5 &1.32e$-$5  \\
\hline
\end{tabular*}
\end{table}

We also conducted simulations using a wide range of bandwidths from
small $h=n^{-1/2}$ to large $h=n^{-1/10}$.
The kernel spectral density curves seem to change rather slowly. This
indicates that the kernel spectral density estimator is robust with
respect to the bandwidth selection.

\section{\texorpdfstring{Proofs of Theorems \protect\ref{theo1} and \protect\ref{theo2}}
{Proofs of Theorems 1 and 2}}

Throughout this section and the next, to simplify notation,
$M,M_1,\ldots,M_{12}$ stand for constants which may take different
values from one appearance to the next.

\subsection{\texorpdfstring{Proof of Theorem \protect\ref{theo1}}{Proof of Theorem 1}}

We begin by developing the following two lemmas, necessary for the
argument of Theorem~\ref{theo1}.
\begin{lemma} \label{lem7}
Under the assumptions of Theorem~\ref{theo1},
 let $F_{c_n,H_n}(t)$ be the distribution function obtained
from $F_{c,H}(t)$
by replacing $c$ and $H$ by $c_n$ and $H_n$, respectively. Furthermore,
$f_{c_n,H_n}(x)$ denotes the density of $F_{c_n,H_n}(x)$. Then,
\[
\sup_{n,x}f_{c_n,H_n}(x)\leq M.
\]
\end{lemma}
\begin{pf}
From (3.10) in~\cite{b4}, we have
%
\begin{equation}\label{a5}
z(\underline{m}_n)=-\frac{1}{\underline{m}_n}+c_n
\int\frac{t\,dH_n(t)}{1+t\underline{m}_n},
\end{equation}
where
$\underline{m}_n=\underline{m}_n(z)=\underline{m}_{F_{c_n,H_n}}(z)$.
Based on this expression, conclusions\vspace*{1pt} similar to those in Theorem 1.1
of~\cite{s1} still hold if we replace $F_{c,H}(x)$ by $F_{c_n,H_n}(x)$
and then argue similarly with the help of~\cite{s1}. For example, the
equality (1.6) in Theorem~1.1 of~\cite{s1} states that
%
\begin{equation}\label{a6}
x=-\frac{1}{\underline{m}(x)}+c\int\frac{t\,dH(t)}{1+t\underline{m}(x)}.
\end{equation}
Similarly, for every $x\neq0$ for which $f_{c_n,H_n}(x)>0$, $\pi
f_{c_n,H_n}(x)$ is the imaginary part of the unique
$\underline{m}_n(x)$ satisfying
%
\begin{equation}\label{a7}
x=-\frac{1}{\underline{m}_n(x)}+c_n\int\frac
{t\,dH_n(t)}{1+t\underline{m}_n(x)}.
\end{equation}

Now, consider the imaginary part of $\underline{m}_n(x).$
From~(\ref{a7}), we obtain
%
\begin{equation}\label{a12}
c_n\int\frac{t^2\,dH_n(t)}{|1+t\underline{m}_n(x)|^2}=\frac{1}{|\underline{m}_n(x)|^2}.
\end{equation}
It follows from~(\ref{a7}),~(\ref{a12}) and H\"older's inequality that
\begin{eqnarray*}
|\underline{m}_n(x)|&\leq&
\frac{|c_n-1|}{x}+\frac{c_n}{x}\int\frac{dH_n(t)}{|1+t\underline{m}_n(x)|}\\
&\leq&\frac{|c_n-1|}{x}+\frac{c_n}{x}\biggl(\int\frac{t^2\,dH_n(t)}
{|1+t\underline{m}_n(x)|^2}\int\frac{dH_n(t)}{t^2}\biggr)^{1/2}\\
&\leq&
\frac{|c_n-1|}{x}+\frac{\sqrt{c_n}}{x|\underline{m}_n(x)|}
\biggl(\int\frac{dH_n(t)}{t^2}\biggr)^{1/2},
\end{eqnarray*}
where $\int\frac{dH_n(t)}{t^2}$ is well defined because we require the
support of $F_{c,H}(x)$ to be $[a,b]$ with $a>0$. This inequality is
equivalent to
\[
|\underline{m}_n(x)|^2\leq\frac{|c_n-1|}{x}|\underline
{m}_n(x)|+\frac{\sqrt{c_n}}{x}\biggl(\int\frac{dH_n(t)}{t^2}\biggr)^{1/2}.
\]
It follows that
%
\begin{equation}\label{a14}
\sup_{n,x}|\underline{m}_n(x)|\leq M.
\end{equation}
This leads to $\sup_{n,x}f_{c_n,H_n}(x)\leq M$.
\end{pf}

\begin{lemma} \label{lem8} Under the assumptions of Lemma
\ref{lem7}, when $x_n\rightarrow x$, we have
%
\begin{equation}\label{a17}
 f_{c_n,H_n}(x_n)-f_{c,H}(x_n)\rightarrow0.
\end{equation}
\end{lemma}

\begin{pf}
Obviously, $f_{c,H}(x_n)-f_{c,H}(x)\rightarrow0$ because $f_{c,H}(x)$
is continuous on the interval $[a,b]$. Moreover, in view of
(\ref{a14}), we may choose a subsequence $n_k$ so that
$\underline{m}_{n_k}(x_{n_k})$ converges. We denote its limit by
$a(x)$. Suppose that $\Im(a(x))>0$. Then, as in Lemma 3.3 in~\cite{s1},
we may argue that the limit of $\underline{m}_{n}(x_{n})$ exists as
$n\rightarrow\infty$. Next, we verify that $a(x)=\underline{m}(x)$. By
(\ref{a7}), we then have
\[
x=-\frac{1}{a(x)}+c\int\frac{t\,dH(t)}{1+ta(x)}
\]
because, via~(\ref{a12}) and H\"older's inequality,
\begin{eqnarray*}
&&\biggl|\int\frac{t\,dH_n(t)}{1+t\underline{m}_n(x)}-\int\frac
{t\,dH_n(t)}{1+ta(x)}\biggr|\\
&&\qquad\leq
|\underline{m}_n(x)-a(x)|\biggl(\frac{1}{c_n|\underline{m}_n(x)|^2}\int
\frac{t^2\,dH_n(t)}{|1+ta(x)|^2}\biggr)^{1/2}
\end{eqnarray*}
and
\[
\int\frac{t\,dH_n(t)}{1+ta(x)}\rightarrow\int\frac{t\,dH(t)}{1+ta(x)}.
\]
Since the solution satisfying the equation~(\ref{a6}) is
unique, $a(x)=\underline{m}(x)$. Therefore,
$\underline{m}_n(x)\rightarrow\underline{m}(x),$ which then implies
that
%
\begin{equation}\label{b37}f_{c_n,H_n}(x_n)-f_{c,H}(x)\rightarrow0.
\end{equation}
Now, suppose that $\Im(a(x))=0$. This implies that
$\Im(\underline{m}_{n}(x_{n}))\rightarrow0$ and then that
$f_{c_n,H_n}(x_n)\rightarrow0$ because if there is another subsequence
on which\break $\Im(\underline{m}_{n}(x_{n}))$ converges to a positive
number, then $\underline{m}_{n}(x_{n})$ must converge to the complex
number with the positive imaginary part, by the previous argument.
Next, by~(\ref{a3}) and~(\ref{a5}),
$\Im(\underline{m}_n(x_n+iv))-\Im(\underline{m}(x_n+iv))\rightarrow0$
for any $v>0$. We may then choose $v_n\rightarrow0$ so that
$\Im(\underline{m}_n(x_n+iv_n))-
\Im(\underline{m}(x_n+iv_n))\rightarrow0$ as $n\rightarrow\infty$.
Moreover, $\Im(\underline{m}(x_n+iv_n))\rightarrow
\Im(\underline{m}(x))$ and $\Im(\underline{m}_n(x_n+iv_n))-
\Im(\underline{m}_n(x))\rightarrow0$ by Theorem~1.1 of~\cite{s1}
and a
theorem for $\underline{m}_n(z)$ similar to Theorem 1.1 of~\cite{s1}.
Therefore, in view of the continuity of $\underline{m}_n(x)$ for
$x\neq
0$, $\Im(\underline{m}(x))=0$ and then~(\ref{a17}) holds for the case
$\Im(a(x))=0$.
\end{pf}

We now proceed to prove Theorem~\ref{theo1}. First, we claim that
%
\begin{equation}\label{a9}
\sup_{x}\biggl|f_n(x)-\frac{1}{h}\int
K\biggl(\frac{x-t}{h}\biggr)\,dF_{c_n,H_n}(t)\biggr|\longrightarrow0
\end{equation}
in probability. Indeed, from integration by parts and Theorem
\ref{theo3}, we obtain
\begin{eqnarray*}
&&E\sup_{x}\biggl|\frac{1}{h}\int
K\biggl(\frac{x-t}{h}\biggr)\,dF^{\bbA_n}(t)-\frac{1}{h}\int
K\biggl(\frac{x-t}{h}\biggr)\,dF_{c_n,H_n}(t)\biggr|\\
&&\qquad=E\sup_{x}\biggl|\frac{1}{h^2}\int
K'\biggl(\frac{x-t}{h}\biggr)\bigl(F^{\bbA_n}(t)-F_{c_n,H_n}(t)\bigr)\,dt\biggr|\\
&&\qquad=E\sup_{x}\biggl|\frac{1}{h}\int
K'(u)\bigl(F^{\bbA_n}(x-uh)-F_{c_n,H_n}(x-uh)\bigr)\,du\biggr|\\
&&\qquad\leq\frac{1}{h}E\sup_{x}|F^{\bbA_n}(x)-F_{c_n,H_n}(x)|\int
|K'(u)|\,du\\
&&\qquad\leq\frac{M}{n^{2/5}h}\rightarrow0.
\end{eqnarray*}
The next aim is to show that
\[
\frac{1}{h}\int
K\biggl(\frac{x-t}{h}\biggr)\,dF_{c_n,H_n}(t)-\frac{1}{h}\int
K\biggl(\frac{x-t}{h}\biggr)\,dF_{c,H}(t) \longrightarrow0
\]
uniformly in $x\in[a,b]$. This is equivalent to, for any sequence
$\{x_n, n\geq1\}$ in $[a,b]$ converging to $x$,
%
\begin{equation}\label{a-10}
\int K(u)\bigl(f_{c_n,H_n}(x_n-uh)-f_{c,H}(x_n-uh)\bigr)\,du
\longrightarrow0.
\end{equation}
From Theorem 1.1 of~\cite{s1},
$f_{c,H}(x)$ is uniformly bounded on the interval $[a,b]$. Therefore,
(\ref{a-10}) follows from the dominated convergence theorem, Lemma~\ref{lem7} and Lemma~\ref{lem8}.

Finally,
\begin{eqnarray*}
&&\biggl|\frac{1}{h}\int K\biggl(\frac{x-t}{h}\biggr)\,dF_{c,H}(t)-f_{c,H}(x)\frac
{1}{h}\int^{x-a}_{x-b}K\biggl(\frac{t}{h}\biggr)\,dt\biggr|\\
&&\qquad=\biggl|\int^{x-a}_{x-b}
\bigl(f_{c,H}(x-t)-f_{c,H}(x)\bigr)\frac{1}{h}K\biggl(\frac{t}{h}\biggr)\,dt\biggr|\\
&&\qquad\leq\sup_{x\in[a,b]}\int_{|t|>\delta} \biggl|\bigl(f_{c,H}(x-t)-f_{c,H}(x)\bigr)
\frac{1}{h}K\biggl(\frac{t}{h}\biggr)\biggr|\,dt\\
&&\qquad\quad{}+\sup_{x\in[a,b]}\int_{|t|\leq\delta}
\biggl|\bigl(f_{c,H}(x-t)-f_{c,H}(x)\bigr)\frac{1}{h}K\biggl(\frac{t}{h}\biggr)\biggr|\,dt\\
&&\qquad\leq2\sup_{x\in[a,b]}f_{c,H}(x)\int_{|t|>
\delta/h}|K(y)|\,dy\\
&&\qquad\quad{}+\sup_{x\in[a,b]}\sup_{|t|\leq\delta}|f_{c,H}(x-t)-f_{c,H}(x)|\int\frac{1}{h}\biggl|K\biggl(\frac{t}{h}\biggr)\biggr|\,dt,
\end{eqnarray*}
which goes to zero by fixing $\delta$ and letting $n\rightarrow\infty$
first, and then letting $\delta\rightarrow0$. On the other hand,
obviously,
\[
\frac{1}{h}\int^{x-a}_{x-b}K\biggl(\frac{t}{h}\biggr)\,dt=\int^{
(x-a)/h}_{(x-b)/h}K(t)\,dt\rightarrow
\int^{+\infty}_{-\infty}K(t)\,dt=1.
\]
Thus, the proof is complete.

\subsection{\texorpdfstring{Proof of Theorem \protect\ref{theo2}}{Proof of Theorem 2}}
Denote by $F_{c_n}(t)$ the distribution function obtained from
$F_{c}(t)=\int_{-\infty}^t f_c(x)\,dx$ with $c$ replaced by $c_n$.
Let $\bbS_n=\frac{1}{n}\bbX_n\bbX_n^T.$ From integration by parts, we
obtain
\begin{eqnarray*}
&&\biggl|\frac{1}{h}\int K\biggl(\frac{x-t}{h}\biggr)\,dF^{\bbS_n}(t)-\frac{1}{h}\int
K\biggl(\frac{x-t}{h}\biggr)\,dF_{c_n}(t)\biggr|\\
&&\qquad=\biggl|\frac{1}{h^2}\int
K'\biggl(\frac{x-t}{h}\biggr)\bigl(F^{\bbS_n}(t)-F_{c_n}(t)\bigr)\,dt\biggr|\\
&&\qquad=\biggl|\frac{1}{h}\int
K'(u)\bigl(F^{\bbS_n}(x-uh)-F_{c_n}(x-uh)\bigr)\,du\biggr|\\
&&\qquad\leq\frac{1}{h}\sup_{x}|F^{\bbS_n}(x)-F_{c_n}(x)|\int|K'(u)|\,du\\
&&\qquad\leq\frac{M}{\sqrt{n}h},
\end{eqnarray*}
where the last step uses Theorem 1.2 in~\cite{g2}. We next prove that
\[
\sup_{x}\biggl|\frac{1}{h}\int
K\biggl(\frac{x-t}{h}\biggr)\,dF_{c_n}(t)-\frac{1}{h}\int
K\biggl(\frac{x-t}{h}\biggr)\,dF_{c}(t)\biggr|\rightarrow0.
\]
It suffices to prove that
%
\begin{equation}\label{a-15}
\sup_{x}|f_{c_n}(x)-f_{c}(x)|\rightarrow0,
\end{equation}
where $f_{c_n}(x)$ stands for the density of $F_{c_n}(x)$.

Note that when $c<1$,
\[
f_{c_n}(x)-f_{c}(x)=\frac{\sqrt{(x-a(c_n))(b(c_n)-x)}}{2\pi
c_nx}-\frac{\sqrt{(x-a(c))(b(c)-x)}}{2\pi cx},
\]
where
\[
a(c)=\bigl(1-\sqrt{c}\bigr)^2,\qquad b(c)=\bigl(1+\sqrt{c}\bigr)^2,
\]
and $a(c_n)$ and $b(c_n)$ are obtained from $a(c)$ and $b(c)$ by
replacing $c$ with $c_n$, respectively.
It is then a simple matter to verify that~(\ref{a-15}) holds for
$x\in[a(c),b(c)]$.

Finally, as in Theorem~\ref{theo1}, one may prove that
\[
\sup_{x}\biggl|\frac{1}{h}\int
K\biggl(\frac{x-t}{h}\biggr)\,dF_{c}(t)-f_{c}(x)\biggr|\rightarrow0.
\]
Thus, the proof is complete.

\subsection{\texorpdfstring{Proof of Corollary \protect\ref{cor1}}{Proof of Corollary 1}}
The result follows from Theorem 1 in~\cite{sch}.

\section{\texorpdfstring{Proof of Theorem \protect\ref{theo3}}{Proof of Theorem 3}}

\subsection{Summary of argument}
 The strategy is to use Corollary 2.2 and
Le\-mma~7.1 in~\cite{g2}. To this end, a key step is to establish an
upper bound for $|b_1|,$ defined below. Note that in a suitable
interval for $z$ with a well-chosen imaginary part $v$, the absolute
value of the expectation of the Stieltjes transform of $F_{\bbA_n}$,
$|Em_n(z)|$, is bounded. Moreover, for such $v$, when
$n\rightarrow\infty$, the difference between $b_1$ and its alternative
expression involving $Em_n(z)$, $\rho_n$ [given in~(\ref{b10})],
converges to zero with some convergence rate. Therefore, we may argue
that $|b_1|$ is bounded. Once this is done, we further develop a
convergence rate of $m_n(z)-Em_n(z)$ using a martingale decomposition,
and a convergence rate of the difference between $Em_n(z)$ and its
corresponding limit using a recurrence approach.

We begin by giving some notation. Define\vspace*{-3pt} $\bbA(z)=\bbA_n-z\bbI$,
$\bbA_j(z)=\bbA(z)-\mathbf{s}_j\mathbf{s}_j^T$ and $\mathbf{s}_j=\bbT
_n^{1/2}\mathbf{x}_j$,
with $\mathbf{x}_j$ being the $j$th column\vadjust{\goodbreak} of $\bbX_n$. Let
$E_j=E(\cdot|\mathbf{s}_1,\ldots,\mathbf{s}_j)$ and let $E_0$ denote the
expectation. Moreover, introduce
\begin{eqnarray*}
\beta_j&=&\frac{1}{1+\mathbf{s}_j^T\bbA_j^{-1}(z)\mathbf{s}_j},\qquad
\hat{\beta}_j=\frac{1}{1+n^{-1}\operatorname{tr}\bbT_n\bbA_j^{-1}(z)},\\
\eta_j&=&\mathbf{s}_j^T\bbA_j^{-1}(z)\mathbf{s}_j-\frac{1}{n}\operatorname{tr}\bbA
_j^{-1}(z)\bbT_n,\qquad
b_1=\frac{1}{1+n^{-1}E\operatorname{tr}\bbT_n\bbA_1^{-1}(z)},\\
m_n(z)&=&\int\frac{dF_{\bbA_n}(x)}{x-z},\qquad
m_n^0(z)=\int\frac{dF_{c_n,H_n}(x)}{x-z},\\
\underline{m}_n(z)&=&\int\frac{dF_{\bbB_n}(x)}{x-z},\qquad
\underline{m}_n^0(z)=\int\frac{d\underline{F}_{c_n,H_n}(x)}{x-z}\\
\end{eqnarray*}
and
\[
\xi_1=\mathbf{s}_1^T\bbA_1^{-1}(z)\mathbf{s}_1-\frac{1}{n}E\operatorname{tr}\bbA
_1^{-1}(z)\bbT_n.
\]
Here, $\underline{F}_{c_n,H_n}(x)$ is obtained from
$\underline{F}_{c,H}(x)$ by replacing $c$ and $H$ by $c_n$ and $H_n$,
respectively.

Let $\Delta_n=\sup_{x}|EF^{\bbA_n}(x)-F_{c_n,H_n}(x)|$ and
$v_0=\max\{\gamma\Delta_n,M_1n^{-2/5} \}$ with $0<\gamma<1$ to be
chosen later and $M_1$ an appropriate constant. As in Lemma~3.1 and
Lemma 3.2 in~\cite{s1}, we obtain, for $u\in[a,b]$ and $v_0\leq v\leq
1,$
%
\begin{equation}\label{b1}
|\underline{m}_n^0(z)|\leq M,\qquad |m_n^0(z)|\leq M,
\end{equation}
where the bound for $|m_n^0(z)|$ is obtained with the help of
(\ref{b3}). Using integration by parts, we have, for $v>v_0$,
\begin{eqnarray*}
|Em_n(z)-m_n^0(z)|&=&\biggl|\int^{+\infty}_{-\infty}\frac{1}{x-z}\,d\bigl(EF^{\bbA
_n}(x)-F_{c_n,H_n}(x)\bigr)\biggr|\\
&=&\biggl|\int^{+\infty}_{-\infty}\frac{EF^{\bbA
_n}(x)-F_{c_n,H_n}(x)}{(x-z)^2}\,dx\biggr|\leq\frac{\pi\Delta_n}{v}\leq
\frac{\pi}{\gamma}.
\end{eqnarray*}
This implies that
%
\begin{equation}\label{b2}
|Em_n(z)|\leq M, \qquad |E\underline{m}_n(z)|\leq M,
\end{equation}
where the bound for $|E\underline{m}_n(z)|$ is obtained from an
equality similar to~(\ref{b3}), noting that $\Re z \geq a$. It is
readily observed that $|\hat{\beta}_j|$ and $|\beta_j|$ are both
bounded by $|z|/v$ (see (3.4) in~\cite{b4}) and that Lemma 2.10 in
\cite{b4} yields
%
\begin{equation}\label{a38}
|\beta_j\mathbf{s}_j^T\bbA_j^{-2}(z)\mathbf{s}_j|\leq v^{-1},
\end{equation}
which gives
%
\begin{equation}\label{b36}
|\operatorname{tr}(\bbA-z\bbI)^{-1}-\operatorname{tr}(\bbA_k-z\bbI)^{-1}|\leq v^{-1}.
\end{equation}
This, together with~(\ref{b2}), gives, for
$v>v_0$,
%
\begin{equation}\label{b9}
\biggl|\frac{1}{n}E\operatorname{tr}\bbA_1^{-1}(z)\biggr|\leq M.\vadjust{\goodbreak}
\end{equation}
In the subsequent subsections, we will assume that $z=u+iv$ with $v\geq
v_0$ and $u\in[a,b]$.

\subsection{\texorpdfstring{Bounds for $n^{-2}E|\operatorname{tr}\bbA^{-1}(z)-E\operatorname{tr}\bbA^{-1}(z)|^2$ and $E|\beta_1|^2$}
{Bounds for $n^{-2}E|\operatorname{tr}\bbA^{-1}(z)-E\operatorname{tr}\bbA^{-1}(z)|^2$ and E|beta1|2}}\label{sub1}

\begin{lemma}\label{lem4}
 If $|b_1|\leq M$, then, for $v>M_1n^{-2/5}$,
%
\begin{equation}\label{a39}
\frac{1}{n^2}E|\operatorname{tr}\bbA^{-1}(z)-E\operatorname{tr}\bbA^{-1}(z)|^2 \leq\frac{M}{n^2v^3}.
\end{equation}
\end{lemma}

\begin{pf}
\begin{eqnarray*}
&&\frac{1}{n}\operatorname{tr}\bbA^{-1}(z)-E\operatorname{tr}\bbA^{-1}(z)\\
&&\qquad=\frac{1}{n}\sum_{j=1}^n\bigl(E_j\operatorname{tr}\bbA^{-1}(z)-E_{j-1}\operatorname{tr}\bbA^{-1}(z)\bigr)\\
&&\qquad=\frac{1}{n}\sum_{j=1}^nE_j\bigl(\operatorname{tr}\bbA^{-1}(z)-\bbA
_j^{-1}(z)\bigr)-E_{j-1}\operatorname{tr}\bigl(\operatorname{tr}\bbA^{-1}(z)-\bbA_j^{-1}(z)\bigr)\\
&&\qquad=\frac{1}{n}\sum_{j=1}^n(E_j-E_{j-1})(\beta_j\mathbf{s}_j^T\bbA
_j^{-2}(z)\mathbf{s}_j)\\
&&\qquad=\frac{1}{n}\sum_{j=1}^n(E_j-E_{j-1})\biggl[b_1\biggl(\mathbf{s}_j^T\bbA
_j^{-2}(z)\mathbf{s}_j-\frac{1}{n}\operatorname{tr}\bbA_j^{-2}(z)\bbT_n\biggr)\\
&&\qquad\quad\hspace*{145pt}{}+b_1\beta
_j\mathbf{s}_j^T\bbA_j^{-2}(z)\mathbf{s}_j\xi_j\biggr],
\end{eqnarray*}
where the last step uses the fact that
%
\begin{equation}\label{b11}
\beta_j=b_1-b_1\beta_j\xi_j.
\end{equation}
Lemma 2.7 in~\cite{b4} then gives
\begin{eqnarray*}
&&E\biggl|\frac{1}{n}\sum_{j=1}^n(E_j-E_{j-1})\biggl(\mathbf{s}_j^T\bbA
_j^{-2}(z)\mathbf{s}_j-\frac{1}{n}\operatorname{tr}\bbA_j^{-2}(z)\bbT_n\biggr)\biggr|^2\\
&&\qquad\leq
\frac{M}{n^2}\sum_{j=1}^nE\biggl|\biggl(\mathbf{s}_j^T\bbA_j^{-2}(z){\bf
s}_j-\frac{1}{n}\operatorname{tr}\bbA_j^{-2}(z)\bbT_n\biggr)\biggr|^2\\
&&\qquad\leq\frac{M}{n^2}\sum_{j=1}^nE\frac{1}{n^2}\operatorname{tr}\bbA
_1^{-2}(z)\bbT_n\bbA_1^{-2}(\bar{z})\bbT_n\\
&&\qquad\leq\frac{\lambda_{\max}^2(\bbT_n)}{n^3v^2}E\operatorname{tr}\bbA_1^{-1}(z)\bbA
_1^{-1}(\bar{z})\leq\frac{M}{n^2v^3}
\end{eqnarray*}
because, via~(\ref{b9}),
%
\begin{equation}\label{b8}
\frac{1}{n}E\operatorname{tr}\bbA_1^{-1}(z)\bbA_1^{-1}(\bar{z})=\frac{1}{v}\Im
\biggl(\frac{1}{n}E\operatorname{tr}\bbA_1^{-1}(z)\biggr)\leq\frac{M}{v}.
\end{equation}
Using~(\ref{a38}) and Lemma 2.7 in~\cite{b4}, we similarly have
\begin{eqnarray*}
&&E\biggl|\frac{1}{n}\sum_{j=1}^n(E_j-E_{j-1})\beta_j{\mathbf
s}_j^T\bbA_j^{-2}(z)\mathbf{s}_j\xi_j\biggr|^2\\
&&\qquad\leq\frac{M}{n^2v^3}+\frac{M}{n^3v^2}E|\operatorname{tr}\bbA^{-1}(z)-E\operatorname{tr}\bbA^{-1}(z)|^2.
\end{eqnarray*}
Summarizing the above, we have proven that
\[
\biggl(1-\frac{M}{nv^2}\biggr)\frac{1}{n^2}E|\operatorname{tr}\bbA^{-1}(z)-E\operatorname{tr}\bbA
^{-1}(z)|^2\leq
\frac{M}{n^2v^3},
\]
which implies Lemma~\ref{lem4} by choosing an appropriate $M_1$ such
that $\frac{M}{nv^2}<\frac{1}{2}$.
\end{pf}

\begin{lemma}\label{lem5}
 If $|b_1|\leq M$, then, for $v>M_1n^{-2/5}$,
%
\begin{equation}\label{a39}
\frac{1}{n^4}E|\operatorname{tr}\bbA^{-1}(z)-E\operatorname{tr}\bbA^{-1}(z)|^4 \leq\frac{M}{n^4v^6}.
\end{equation}
\end{lemma}
\begin{pf}
Lemma~\ref{lem5} is obtained by repeating the argument of Lemma
\ref{lem4} and applying
\begin{eqnarray*}
&&E\biggl(\frac{1}{n}\operatorname{tr}\bbA_1^{-2}(z)\bbT_n\bbA_1^{-2}(\bar{z})\bbT
_n\biggr)^2\\
&&\qquad\leq
\frac{\lambda_{\max}^4(\bbT_n)}{n^2v^6}E|\operatorname{tr}\bbA_1^{-1}(z)-E\operatorname{tr}\bbA
_1^{-1}(z)|^2+\frac{\lambda_{\max}^4(\bbT_n)}{n^2v^6}|E\operatorname{tr}\bbA_1^{-1}(z)|^2\\
&&\qquad\leq
\frac{M}{v^6}.
\end{eqnarray*}\upqed\hspace*{-280pt}
\end{pf}

\begin{lemma}\label{lem6}
If $|b_1|\leq M$, then there is some
constant $M_2$ such that for $v\geq M_2n^{-2/5}$,
\[
E|\beta_1|^2\leq M.
\]
\end{lemma}

\begin{pf}
By~(\ref{b11}), we have
\[
\beta_j=b_1-b_1^2\xi_j+b_1^2\beta_j\xi_j
\]
and
\begin{eqnarray}\label{b12}\qquad
E|\xi_1(z)|^4&\leq& ME|\eta_1(z)|^4+Mn^{-4}E|\operatorname{tr}\bbA_1^{-1}(z)\bbT_n-E\operatorname{tr}\bbA
_1^{-1}(z)\bbT_n|^4\nonumber\\ [-8pt]\\ [-8pt]
&\leq& \frac{M}{n^2v^2}+\frac{M}{n^4v^6}\nonumber
\end{eqnarray}
because repeating the argument of Lemma~\ref{lem4} and Lemma~\ref{lem5}
yields
%
\begin{equation}\label{b38}
E\biggl|\frac{1}{n}\operatorname{tr}\bbD\bbA_1^{-1}(z)-E\frac{1}{n}\operatorname{tr}\bbD\bbA
_1^{-1}(z)\biggr|^4\leq\frac{M}{n^4v^{6}\|\bbD\|^4}
\end{equation}
for a fixed matrix $\bbD$. It follows that
\[
E|\beta_1|^2\leq
|b_1|^2+|b_1|^4E|\xi_1|^2+\frac{|b_1|^4}{v}(E\vert\beta_1|^2E|\xi_1|^4)^{1/2},
\]
which gives
\[
E|\beta_1|^2\leq M+\frac{M}{nv}+\frac{M}{nv^2}(E\vert\beta_1|^2)^{1/2}.
\]
Solving this inequality gives Lemma~\ref{lem6}.\vspace*{-2pt}
\end{pf}

\subsection{A bound for $b_1(z)$}

By~(\ref{b11}) and
%
\begin{equation}\label{b7}
1-c_n-zc_nm_n(z)=\frac{1}{n}\sum_{j=1}^n\beta_j
\end{equation}
(see the equality above (2.2) in~\cite{s3}), we get
%
\begin{equation}\label{b10}
b_1=1-c_n-zc_nEm_n(z)+\rho_n,
\end{equation}
where\vspace*{-2pt}
\[
\rho_n=b_1E(\beta_1\xi_1).\vspace*{-2pt}
\]

\begin{lemma}\label{lem2}
If $|b_1|\leq M$, then there is some
constant $M_3$ such that for $v\geq M_3n^{-2/5}$,
\[
|\rho_n|\leq\frac{M}{nv}.\vspace*{-2pt}
\]
\end{lemma}

\begin{pf}
Lemma~\ref{lem6} and~(\ref{b12}) ensure that
\[
|E[\beta_1(z)\xi_1(z)]|=|b_1(z)E[\beta_1(z)\xi_1^2]| \leq
M(E|\beta_1(z)|^2E|\xi_1|^4)^{1/2}\leq\frac{M}{nv}.
\]
Thus, Lemma~\ref{lem2} is proved.\vspace*{-2pt}
\end{pf}

\begin{lemma}\label{lem1}
If $\Im(z+\rho_n)\geq0,$ then there exists a positive constant $c$
depending on $\gamma,a,b$ such that
\[
|b_1|\leq M.\vspace*{-2pt}
\]
\end{lemma}

\begin{pf}
Consider the case $\Im(Em_n(z))\geq v>0$ first. It follows from~(\ref{b10}) and the assumption that
\begin{eqnarray*}
\Im\bigl(c_n+z+zc_nEm_n(z)-1\bigr)
&\geq&-\Im(b_1)\\
&=&-|b_1|^2\Im\bigl(1+n^{-1}E\operatorname{tr}\bbA^{-1}(\bar{z})\bigr).\vadjust{\goodbreak}
\end{eqnarray*}
Note that
%
\begin{eqnarray}\label{c4}
&&\Im\bigl(c_n+z+zc_nEm_n(z)-1\bigr)\nonumber\\
&&\qquad= v+vc_n\int\frac
{x}{|x-z|^2}\,dF_{n2}(x)\\
&&\qquad=v+c_n[v\Re(Em_n(z))+u\Im(Em_n(z))]>0.\nonumber
\end{eqnarray}
Thus, we have
%
\begin{eqnarray}\label{b33}
|b_1|^2&\leq&\frac{v+c_n[v\Re(Em_n(z))+u\Im(Em_n(z))]}{c_n\Im(Em_n(z))}\nonumber\\
&\leq&\frac{[1+c_n|\Re(Em_n(z))|+c_nu]\Im(Em_n(z))}{c_n\Im(Em_n(z))}\\
&\leq&1/c_n+M+b.\nonumber
\end{eqnarray}

Next, consider the case $\Im(Em_n(z))<v$. Note that for $u\in[a,b]$,
%
\begin{equation}\label{b18}
|\Im(E\underline{m}_n(z))|\geq\frac{v}{M+v^2}.
\end{equation}
This, together with~(\ref{b33}), gives
\[
\hspace*{-15pt}|b_1|^2\leq\frac{(M+v^2)[1+c_n(|\Re(Em_n(z))|+u)]v}{c_nMv}
\leq\frac{1+c_n[|\Re(Em_n(z))|+u]}{c_nM}.
\]
\upqed\end{pf}

\begin{lemma}\label{lem3}
There is some constant $M_4$ such that, for any $ v\geq M_4n^{-2/5}$,
\[
\Im(z+\rho_n)> 0.
\]
\end{lemma}

\begin{pf}
First, we claim that
%
\begin{equation}\label{a10}
\Im(z+\rho_n)\neq0.
\end{equation}
If not, $\Im(z+\rho_n)=0$ implies that
%
\begin{equation}\label{a8}
|\rho_n|\geq|\Im(\rho_n)|=v.
\end{equation}
On the
other hand, if $\Im(z+\rho_n)=0$, then we then conclude from Lemma~\ref{lem1} and Lemma~\ref{lem2} that
\[
|\rho_n|\leq\frac{M}{nv}.
\]
Thus, recalling that $v\geq M_4n^{-2/5}$, we may choose an appropriate
constant $M_4$ so that
\[
|\rho_n|\leq\frac{v}{3},
\]
which contradicts~(\ref{a8}). Therefore,~(\ref{a10}) holds.

Next, note that
\[
\Im\bigl(z+zn^{-1}E\operatorname{tr}\bbA_1^{-1}(z)\bigr)\geq v,\qquad\Im\bigl(z+zn^{-1}E\operatorname{tr}\bbA^{-1}(z)\bigr)\geq v.
\]
Therefore, when taking $v=1$,
\[
|b_1(z)|\leq\frac{|z|}{v}\leq M,\qquad|b(z)|\leq\frac{|z|}{v}\leq M.
\]
It follows from Lemma~\ref{lem1} and Lemma~\ref{lem2} that
\[
|\rho_n|\leq\frac{M}{n},
\]
which implies that for $n$ large and $v=1$,
%
\begin{equation}\label{b13}
\Im(z+\rho_n)> 0.
\end{equation}
This, together with~(\ref{a10}) and continuity of the function, ensures
that~(\ref{b13}) holds for $1\geq v\geq M_3n^{-2/5}$. Thus, the proof
of Lemma~\ref{lem3} is complete.
\end{pf}

\subsection{Convergence of expected value}
Based on Lemma~\ref{lem1} and Lemma~\ref{lem3}, $|b_1|\leq M$ and
therefore all results in Section~\ref{sub1} remain true for $v\geq
Mn^{-2/5}$ with some appropriate positive constant $M$.

Set
$ \bbF^{-1}(z)=(E\underline{m}_n\bbT_n+\bbI)^{-1}$ and then write (see
(5.2) in~\cite{b4})
%
\begin{equation}\label{b15}
c_n\int\frac{dH_n(t)}{1+tE\underline{m}_n}+zc_nE(m_n(z))=D_n,
\end{equation}
where
\[
D_n=E\beta_1\biggl[\mathbf{s}_1^T\bbA_1^{-1}(z)\bbF^{-1}(z){\bf
s}_1-\frac{1}{n}E(\operatorname{tr}\bbF^{-1}(z)\bbT_n\bbA^{-1}(z))\biggr].
\]
It follows that (see (3.20) in~\cite{b4})
\begin{eqnarray}\label{a32}\quad
&&E\underline{m}_n(z)-\underline{m}_n^0(z)\nonumber\\ [-8pt]\\ [-8pt]
&&\qquad=\underline
{m}_n^0(z)E\underline{m}_n\omega_n\Big/\biggl(1-c_nE\underline{m}_n\underline
{m}_n^0\int\frac{t^2\,dH_n(t)}{(1+tE\underline{m}_n)(1+t\underline{m}_n^0)}\biggr),\nonumber
\end{eqnarray}
where $ \omega_n=-D_n/E\underline{m}_n.$

Applying~(\ref{b11}), we obtain
\begin{eqnarray*}
D_n&=&b_1E\biggl[\frac{1}{n}\operatorname{tr}\bbF^{-1}(z)\bbT_n\bbA_1^{-1}(z)-\frac
{1}{n}\operatorname{tr}\bbF^{-1}(z)\bbT_n\bbA^{-1}(z)\biggr]\\
&&{}-E\biggl[b_1\beta_1\xi_1
\biggl(\mathbf{s}_1^T\bbA_1^{-1}(z)\bbF^{-1}(z)\mathbf{s}_1-\frac{1}{n}E(\operatorname{tr}\bbF
^{-1}(z)\bbT_n\bbA^{-1}(z))\biggr)\biggr].
\end{eqnarray*}

We now investigate $D_n$. We conclude from~(\ref{a38}) and H\"older's
inequality that
\begin{eqnarray}\label{a33}
&&\biggl|\frac{1}{n}\operatorname{tr}\bbF^{-1}(z)\bbT_n\bbA_1^{-1}(z)-\frac{1}{n}\operatorname{tr}\bbF
^{-1}(z)\bbT_n\bbA^{-1}(z)\biggr|\nonumber\\ [-8pt]\\ [-8pt]
&&\qquad\leq
\frac{M}{nv^{3/2}}\biggl(\frac{1}{n}\operatorname{tr}\bbF^{-1}(z)\bbF^{-1}(\bar{z})\biggr)^{1/2}.\nonumber
\end{eqnarray}
Let
$\zeta_1=\mathbf{s}_1^T\bbA_1^{-1}(z)\bbF^{-1}(z)\mathbf{s}_1-\frac
{1}{n}(\operatorname{tr}\bbF^{-1}(z)\bbT_n\bbA_1^{-1}(z))$.
By~(\ref{b8}) and H\"older's inequality, we have
\[
|Eb_1^2\xi_1 \zeta_1|=|b_1^2E\eta_1\zeta_1|\leq
\frac{M}{nv^{3/2}}\biggl(\frac{1}{n}\operatorname{tr}\bbF^{-1}(z)\bbF^{-1}(\bar{z})\biggr)^{1/2}
\]
and by Lemma~\ref{lem6}, Lemma~\ref{lem5},~(\ref{b12}),~(\ref{b14}) and
H\"older's inequality, we have
\begin{eqnarray*}
&&E|b_1^2\beta_1\xi_1^2 \zeta_1|\\
&&\qquad\leq M
(E|\beta_1|^2)^{1/2}(E|\eta_1|^8E|\zeta_1|^4)^{1/4}\\
&&\qquad\quad{}+M(E|\beta_1|^2)^{1/2}\\
&&\qquad\quad\quad{}\times\biggl(E\biggl[\biggl|\frac{1}{n}\operatorname{tr}\bbA_1^{-1}(z)\bbT_n-E\frac
{1}{n}\operatorname{tr}\bbA_1^{-1}(z)\bbT_n\biggr|^4E(|\zeta_1|^2|\bbA
_1^{-1}(z))\biggr]\biggr)^{1/2}\\
&&\qquad\leq\frac{M}{nv^{3/2}},
\end{eqnarray*}
where we also use~(\ref{b38}) and the fact that, via Lemma 2.11 in
\cite{b4},
%
\begin{equation}\label{b14}
\|\bbF^{-1}(z)\|\leq\frac{M}{v}.
\end{equation}
These, together with~(\ref{b11}),
give
\begin{eqnarray}\label{a34}
|Eb_1\beta_1\xi_1 \zeta_1|&\leq&|Eb_1^2\xi_1 \zeta_1|+
|Eb_1^2\beta_1\xi_1^2 \zeta_1|\nonumber\\ [-8pt]\\ [-8pt]
&\leq&\frac{M}{nv^{3/2}}+\frac{M}{nv^{3/2}}\biggl(\frac{1}{n}\operatorname{tr}\bbF
^{-1}(z)\bbF^{-1}(\bar{z})\biggr)^{1/2}.\nonumber
\end{eqnarray}
Similarly, by~(\ref{b38}), we may get
%
\begin{equation}\label{a36}
 \biggl|Eb_1\beta_1\xi_1
\biggl(\frac{1}{n}\operatorname{tr}\bbF^{-1}(z)\bbT_n\bbA_1^{-1}(z)\!-\!E\frac{1}{n}\operatorname{tr}\bbF
^{-1}(z)\bbT_n\bbA_1^{-1}(z)\biggr)\biggr|
\leq\frac{M}{nv^{3/2}}.\hspace*{-40pt}
\end{equation}
In view of~(\ref{a33}), we have
\begin{eqnarray*}
&&E\biggl|b_1\beta_1\xi_1 \biggl(E
\frac{1}{n}\operatorname{tr}\bbF^{-1}(z)\bbT_n\bbA_1^{-1}(z)-E\frac{1}{n}\operatorname{tr}\bbF
^{-1}(z)\bbT_n\bbA^{-1}(z)\biggr)\biggr|\\
&&\qquad\leq
\frac{M}{nv}\biggl(\frac{1}{n}\operatorname{tr}\bbF^{-1}(z)\bbF^{-1}(\bar{z})\biggr)^{1/2}.
\end{eqnarray*}
Summarizing the above gives
%
\begin{equation}\label{b16}
|D_n|\leq\frac{M}{nv^{3/2}}+\frac{M}{nv^{3/2}}\biggl(\frac{1}{n}\operatorname{tr}\bbF
^{-1}(z)\bbF^{-1}(\bar{z})\biggr)^{1/2}.
\end{equation}

Now, considering the imaginary part of~(\ref{b15}), we may conclude
that
%
\begin{equation}\label{b17}
c_n\int\frac{t}{|1+tE\underline{m}_n|^2}\,dH_n(t)\leq\frac{|\Im
(zc_nE(m_n(z)))|}{\Im(E\underline{m}_n)}+\frac{|D_n|}{\Im
(E\underline{m}_n)}.
\end{equation}
Formulas~(\ref{b18}),~(\ref{b2}) and an equality similar to~(\ref{b3})
ensure that
%
\begin{equation}\label{b19}
 \frac{|\Im(zc_nE(m_n(z)))|}{\Im(E\underline{m}_n)}\leq
\frac{u\Im(E\underline{m}_n)+v|\Re(E\underline{m}_n)|}{\Im
(E\underline{m}_n)}\leq
M
\end{equation}
and that
\[
\frac{|D_n|}{\Im(E\underline{m}_n)}\leq
\frac{M}{nv^{5/2}}+\frac{M}{nv^{5/2}}\biggl(\frac{1}{n}\operatorname{tr}\bbF^{-1}(z)\bbF
^{-1}(\bar{z})\biggr)^{1/2}.
\]
It follows that
\begin{eqnarray}\label{b20}
&&c_n\int\frac{t}{|1+tE\underline{m}_n|^2}\,dH_n(t)\nonumber\\ [-8pt]\\ [-8pt]
&&\qquad\leq
M+\frac{M}{nv^{5/2}}+\frac{M}{nv^{5/2}}\biggl(\frac{1}{n}\operatorname{tr}\bbF
^{-1}(z)\bbF^{-1}(\bar{z})\biggr)^{1/2},\nonumber
\end{eqnarray}
which implies that
\begin{eqnarray*}
\biggl|\frac{1}{n}\operatorname{tr}\bbF^{-1}(z)\bbF^{-1}(\bar{z})\biggr|&=&\int\frac
{dH_n(t)}{|1+tE\underline{m}_n|^2}\leq
\frac{1}{\lambda_{\min}(\bbT_n)}\int\frac
{t\,dH_n(t)}{|1+tE\underline{m}_n|^2}\\
&\leq& M+\frac{M}{nv^{5/2}}+\frac{M}{nv^{5/2}}\biggl(\frac{1}{n}\operatorname{tr}\bbF
^{-1}(z)\bbF^{-1}(\bar{z})\biggr)^{1/2}.
\end{eqnarray*}
This inequality yields
%
\begin{equation}\label{b21}
\biggl|\frac{1}{n}\operatorname{tr}\bbF^{-1}(z)\bbF^{-1}(\bar{z})\biggr|\leq M.
\end{equation}
This, together with~(\ref{b16}), ensures that
%
\begin{equation}\label{b22}
|D_n|\leq\frac{M}{nv^{3/2}}.
\end{equation}

Next, we prove that
%
\begin{equation}\label{b23}
\inf_{n,z} |E\underline{m}_n(z)|>M>0.
\end{equation}
To this end, by~(\ref{b10}) and an equality similar to~(\ref{b3}), we
have
%
\begin{equation}
b_1=-zE\underline{m}_n(z)+\rho_n.\label{b25}
\end{equation}
In view of Lemma~\ref{lem2} and~(\ref{b25}), to prove~(\ref{b23}),
it is
thus sufficient to show that
%
\begin{equation}\label{b24}
\biggl|\frac{1}{n}E\operatorname{tr}\bbA_1^{-1} (z)\bbT_n\biggr|\leq M.
\end{equation}
Suppose that~(\ref{b24}) is not true. There then exist subsequences
$n_k$ and $z_k\rightarrow z_0\neq0$ such that
$|\frac{1}{n}E\operatorname{tr}\bbA_1^{-1} (z)\bbT_{n}|\rightarrow\infty$ on the
subsequences $n_k$ and $z_k$, which, together with~(\ref{b25}) and
Lemma~\ref{lem2}, implies that $E\underline{m}_n(z)\rightarrow0$ on
such subsequences. This, together with~(\ref{b21}), ensures that on
such subsequences
\[
c_n\int\frac{dH_n(t)}{1+tE\underline{m}_n}\rightarrow c,
\]
which, via an equality similar to~(\ref{b3}), further implies that on
such subsequences,
%
\begin{equation}\label{b26}
c_n\int\frac{dH_n(t)}{1+tE\underline{m}_n}+zc_nE(m_n(z))\rightarrow1.
\end{equation}
But, on the other hand, by~(\ref{b22}) and~(\ref{b15}),
\[
c_n\int\frac{dH_n(t)}{1+tE\underline{m}_n}+zc_nE(m_n(z))\rightarrow0,
\]
which contradicts~(\ref{b26}). Therefore,\vspace*{1pt}~(\ref{b24}) and,
consequently,~(\ref{b23}) hold.

It follows from~(\ref{b23}) and~(\ref{b22}) that for $v>M_8n^{-2/5}$,
%
\begin{equation}\label{b27}
|\omega_n|\leq\frac{M}{nv^{3/2}}\leq v,
\end{equation}
where we may choose an appropriate $M_8$. Moreover, since~(\ref{a4})
holds when $\underline{m}$ is replaced by $\underline{m}_n^0$,
considering the imaginary parts of both sides of the equality, we
obtain
\[
v=\frac{\Im(\underline{m}_n^0)}{|\underline{m}_n^0|^2}-c_n\Im
(\underline{m}_n^0)\int\frac{t^2\,dH_n(t)}{|1+t\underline{m}_n^0|^2},
\]
which implies that
\[
c_n\Im(\underline{m}_n^0)\int\frac{t^2\,dH_n(t)}{|1+t\underline
{m}_n^0|^2}\leq
M.
\]
It follows that
\[
\biggl(\biggl(c_n\Im(\underline{m}_n^0)\int\frac
{t^2\,dH_n(t)}{|1+t\underline{m}_n^0|^2}\biggr)\Big/\biggl(v+c_n\Im(\underline
{m}_n^0)\int\frac{t^2\,dH_n(t)}{|1+t\underline{m}_n^0|^2}\biggr)\biggr)^{1/2}
\leq1-Mv.
\]
Applying this and~(\ref{b27}), as in (3.21) in~\cite{b4}, we may
conclude that
%
\begin{equation}\label{b28}
\biggl|1-c_nE\underline{m}_n\underline{m}_n^0\int\frac
{t^2\,dH_n(t)}{(1+tE\underline{m}_n)(1+t\underline{m}_n^0)}\biggr|\geq
Mv.
\end{equation}
This, together with~(\ref{a32}) and~(\ref{b22}), yields
%
\begin{equation}\label{b29}
|E\underline{m}_n(z)-\underline{m}_n^0(z)|\leq\frac{M}{nv^{5/2}}.
\end{equation}

\subsection{Convergence rate of $EF^{\bbA_n}$ and $F^{\bbA_n}$}

As in Theorem 1.1 in~\cite{s1},
$f_{c_n,H_n}$ is continuous. Therefore,
\[
\frac{1}{v\pi} \sup_{x\in
[a+1/2\varepsilon,b-1/2\varepsilon]}\int
_{|y|<2vM}|F_{c_n,H_n}(x+y)-F_{c_n,H_n}(x)|\,dy\leq
Mv,
\]
where $\varepsilon> vM_{11}$. Lemma 2.1 in~\cite{g1} or Lemma 2.1 and
Corollary 2.2 in~\cite{g2} are then applicable in our case.

First, consider $EF^{\bbA_n}$. For $v\geq v_0$, by Corollary 2.2 in
\cite{g2},~(\ref{b29}), we obtain, after integration in $u$ and $v,$
%
\begin{equation}\label{b30}
\Delta_n\leq\frac{M}{n}+M_9v_0+\frac{M_{10}}{nv_0^{3/2}},
\end{equation}
where we set $V$, given in Corollary 2.2 in~\cite{g2}, equal to one and
also use the fact that
$|E\underline{m}_n(z')-\underline{m}_n^0(z')|=O(n^{-1})$ with $z'=u+iV$
(see Section 4 in~\cite{b1}). If $v_0=M_1n^{-2/5}$, then~(\ref{b30})
gives $|\Delta_n|\leq M/n^{2/5}$. If $v_0=\gamma\Delta_n$, then we
choose $\gamma=(2M_9)^{-1}$ (here one should note that $M_{10}$ depends
on $\gamma$, but $M_9$ does not depend on $\gamma$). Again,~(\ref{b30})
gives
%
\begin{equation}\label{b31}
|\Delta_n|\leq M/n^{2/5}.
\end{equation}
This completes the proof of
(\ref{b34}).

Now, consider the convergence rate of $F^{\bbA_n}$. It follows from
Cauchy's inequality that
\[
n^{-1}|\operatorname{tr}\bbA^{-2}(z)-E\operatorname{tr}\bbA^{-2}(z)|\leq\frac{M}{v}\sup
_{z_1\in
\mathcal{C}_v}n^{-1}|\operatorname{tr}\bbA^{-1}(z_1)-E\operatorname{tr}\bbA^{-1}(z_1)|,
\]
where $\mathcal{C}_v=\{z_1\dvtx|z-z_1|=v_0/3\}$. This, together with Lemma
\ref{lem4}, ensures that
%
\begin{equation}\label{b32}
En^{-1}|\operatorname{tr}\bbA^{-2}(z)-E\operatorname{tr}\bbA^{-2}(z)|\leq\frac{M}{nv^{5/2}}.
\end{equation}
Equation~(\ref{b35}) then follows from~(\ref{b32}), Lemma~\ref{lem4}, the
argument leading to~(\ref{b31}) and Lemma 7.1 in~\cite{g2}.

\section*{Acknowledgments} The authors would like to thank the Editor, an
Associate Editor and a referee for their constructive comments which
helped to
improve this paper considerably.


\printaddresses


\begin{thebibliography}{16}

\bibitem{b3}
\textsc{Bai}, Z. D. (1993).
Convergence rate of expected spectral
distributions of large random matrices. Part II. Sample covariance
matrices.
\textit{Ann. Probab.} \textbf{21} 649--672.
\MR{1217560}\vadjust{\goodbreak}

\bibitem{b4}
\textsc{Bai}, Z. D. and \textsc{Silverstein}, J. W. (1998).
No eigenvalues outside the
support of the limiting spectral distribution of large dimensional
random matrices.
\textit{Ann. Probab.} \textbf{26} 316--345.
\MR{1617051}

\bibitem{b1}
\textsc{Bai}, Z. D. and \textsc{Silverstein}, J. W. (2004).
CLT for linear spectral
statistics of large dimensional sample covariance matrices.
\textit{Ann. Probab.} \textbf{32} 553--605.
\MR{2040792}

\bibitem{bai06}
\textsc{Bai}, Z. D. and \textsc{Silverstein}, J. W. (2006).
\textit{Spectral Analysis of Large Dimensional Random Matrices.
Mathematics Monograph Series} \textbf{2}. Science Press, Beijing.

\bibitem{g1}
\textsc{G\"otze}, F. and \textsc{Tikhomirov}, A. (2003).
Rate of convergence to the semi-circular law.
\textit{Probab. Theory Related. Fields} \textbf{127} 228--276.
\MR{2013983}

\bibitem{g2}
\textsc{G\"otze}, F. and \textsc{Tikhomirov}, A. (2004).
Rate of convergence in probability to the Mar\v{c}enko--Pastur law.
\textit{Bernoulli} \textbf{10} 503--548.
\MR{2061442}

\bibitem{peter}
\textsc{Hall}, P. (1984).
An optimal property of kernel estimators of a probability
density.
\textit{J. Roy. Statist. Soc. Ser. B} \textbf{1} 134--138.
\MR{0745225}

\bibitem{MP}
\textsc{Mar\v{c}enko}, V. A. and \textsc{Pastur}, L. A. (1967).
Distribution for some sets of random matrices.
\textit{Math. USSR-Sb.} \textbf{1} 457--483.

\bibitem{p1}
\textsc{Parzen}, E. (1962).
On estimation of a probability density function and mode.
\textit{Ann. Math. Stat.} \textbf{33} 1065--1076.
\MR{0143282}

\bibitem{R1}
\textsc{Rosenblatt}, M. (1956).
Remarks on some non-parametric estimates of a density function.
\textit{Ann. Math. Stat.} \textbf{3} 832--837.
\MR{0079873}

\bibitem{sch}
\textsc{Scheff}, H. (1947).
A useful converge theorem for probability distributions.
\textit{Ann. Math. Stat.} \textbf{3} 434--438.
\MR{0021585}

\bibitem{bw}
\textsc{Silverman}, B. W. (1986).
Nonparametric estimation for statistics and data analysis.
Chapman \& Hall, New York.
\MR{0848134}

\bibitem{s3}
\textsc{Silverstein}, J. W. (1995).
Strong convergence of the limiting distribution of the eigenvalues of
large dimensional random
matrices.
\textit{J. Multivariate Anal.} \textbf{55} 331--339.
\MR{1370408}

\bibitem{s1}
\textsc{Silverstein}, J. W. and \textsc{Choi}, S. I. (1995).
Analysis of the limiting spectral distribution of large dimensional
random matrices.
\textit{J. Multivariate Anal.} \textbf{54} 295--309.
\MR{1345541}

\bibitem{Tse99}
\textsc{Tse}, D. and \textsc{Hanly}, S. (1999).
Linear multiuser receivers: Effective interference, effective bandwidth
and user capacity.
\textit{IEEE Trans. Inform. Theory} \textbf{45} 641--657.
\MR{1677023}

\bibitem{y1}
\textsc{Yin}, Y. Q. (1986).
Limiting spectral distribution for a class of random matrices.
\textit{J.~Multivariate Anal.} \textbf{20} 50--68.
\MR{0862241}

\end{thebibliography}
\end{document}